\documentclass[12pt]{amsart}

\usepackage{amsmath}
\usepackage{amssymb}
\usepackage{amsthm}
\usepackage[pdftex]{graphicx}

\theoremstyle{definition}
\newtheorem{thm}{Theorem}
\newtheorem{cor}[thm]{Corollary}
\newtheorem{lemma}[thm]{Lemma}
\newtheorem{prop}[thm]{Proposition}
\newtheorem{remark}[thm]{Remark}
\newtheorem{conjecture}[thm]{Conjecture}
\newtheorem*{defn}{Definition}

\newtheorem*{claimm}{Claim}

\newcommand{\fieldk}{k}
\newtheorem*{EGH}{EGH Conjecture}
\DeclareMathOperator{\soc}{soc}
\DeclareMathOperator{\depth}{\text{depth}}

\begin{document}

\title[An application of liaison theory to the EGH conjecture]{An application of liaison theory to the Eisenbud-Green-Harris conjecture}
\author{Kai Fong Ernest Chong}
\address{Department of Mathematics\\
      Cornell University\\
      Ithaca, NY 14853-4201, USA}
\email{kc343@cornell.edu}


\begin{abstract}
In this paper, we apply liaison theory to the Eisenbud-Green-Harris conjecture and prove that the conjecture holds for a certain subclass of homogeneous ideals in the linkage class of a complete intersection ideal. In the case of three variables, we prove that the conjecture holds for Gorenstein ideals.
\end{abstract}

\maketitle

\section{Introduction}\label{sec:Intro}
The Eisenbud-Green-Harris (EGH) conjecture is a famous open problem in commutative algebra and algebraic geometry. If true, it would imply the generalized Cayley-Bacharach conjecture~\cite{EisenbudGreenHarris1993:EGHConjecture}, generalize the classic Clements-Lindstr\"{o}m theorem~\cite{ClementsLindstrom1969} in combinatorics, as well as extend Macaulay's characterization~\cite{Macaulay1927} of the Hilbert functions of homogeneous ideals in a polynomial ring $S:= \fieldk[x_1, \dots, x_n]$ to the case of homogeneous ideals containing a given $S$-regular sequence.

There are several equivalent formulations and slight variations of the EGH conjecture in the literature; see \cite{FranciscoRichert2007:LPPIdeals} for a good overview. Given a proper ideal $I$ of a Noetherian ring $R$, we say $I$ {\it minimally contains an $(a_1, \dots, a_r)$-regular sequence of forms} (i.e. homogeneous polynomials) if $I$ has depth $r$, the minimum degree of the forms in $I$ is $a_1$, and for each $2\leq i \leq r$, the integer $a_i$ is the smallest degree such that $I$ contains an $R$-regular sequence $f_1, \dots, f_i$ of forms of degrees $a_1, \dots, a_i$ respectively. In this paper, we work with the following version of the conjecture.

\EGH Let $2\leq e_1 \leq \dots \leq e_n$ be integers. If $I\subsetneq S$ is a homogeneous ideal that minimally contains an $(e_1, \dots, e_n)$-regular sequence of forms, then there exists a homogeneous ideal $J\subsetneq S$ containing $x_1^{e_1}, \dots, x_n^{e_n}$, such that $I$ and $J$ have the same Hilbert function.

The EGH conjecture is known to be true in some cases. Richert~\cite{Richert2004:StudyLPPConjecture} settled the case $n=2$, Francisco~\cite{Francisco2004:AlmostCIandLPPConj} proved the case when $I$ is an almost complete intersection ideal, Caviglia and Maclagan~\cite{CavigliaMaclagan2008:EGHConj} showed that the conjecture is true when $e_{j+1} > \sum_{i=1}^j (e_i-1)$ for all $1\leq j<n$, and Abedelfatah~\cite{Abedelfatah2012:EGHConjPreprint} proved the conjecture when every $f_i$ is a product of linear forms. In the special case when $e_i = 2$ for all $1\leq i\leq n$, Richert~\cite{Richert2004:StudyLPPConjecture} claimed in unpublished work that the conjecture holds for $n\leq 5$, Chen~\cite{Chen2011:thesis} gave a proof for $n\leq 4$, Herzog and Popescu~\cite{HerzogPopescu1998:HilbertFunctionsGenericForms} proved that the conjecture holds if $\fieldk$ is a field of characteristic zero and $I$ is minimally generated by generic quadratic forms, while Gasharov~\cite{Gasharov1999:HilbertFunctionsHomogGenericForms} showed that Herzog-Popescu's result is still true when $\fieldk$ is replaced by a field of arbitrary characteristic. As for the case $n=3$, Cooper~\cite{Cooper:EGHPreprintCase3} proved the conjecture when $e_1 \leq 3$ by considering the remaining cases not covered by Caviglia-Maclagan's result.

In contrast to these known results, we use liaison theory as our main tool. Licci ideals are homogeneous ideals in the liaison class of a complete intersection ideal, and projective schemes defined by licci ideals are fundamental objects studied in liaison theory. In this paper, we introduce a new subclass of licci ideals, which we call `sequentially bounded', and we prove that the EGH conjecture holds for every sequentially bounded licci ideal that `admits a minimal first link' (Theorem \ref{thm:SeqBoundedLicci}); see Section \ref{sec:SeqBoundedLicci} for a precise definition of such licci ideals. As an important consequence, we show that the EGH conjecture holds for Gorenstein ideals in the case of three variables:

\thm\label{mainThm-Gorenstein} Let $2\leq e_1 \leq e_2 \leq e_3$ be integers. If $I\subsetneq \fieldk[x_1,x_2,x_3]$ is a homogeneous Gorenstein ideal that minimally contains an $(e_1, e_2, e_3)$-regular sequence of forms, then there exists a monomial ideal $J$ in $\fieldk[x_1,x_2,x_3]$ containing $x_1^{e_1}, x_2^{e_2}, x_3^{e_3}$, such that $I$ and $J$ have the same Hilbert function.

Our proof of Theorem \ref{mainThm-Gorenstein} uses Migliore-Nagel's recent result~\cite{MiglioreNagel2010:MinimalLinksResultGaeta} in liaison theory, which in turn relies on the Buchsbaum-Eisenbud structure theorem~\cite{BuchsbaumEisenbud1977} on Gorenstein ideals of height three. We hope that further advances in liaison theory would give more insight into proving the EGH conjecture.

\section{Liaison Theory}\label{sec:LiaisonTheory}
Let $\mathbb{N}$ and $\mathbb{P}$ denote the non-negative and positive integers respectively. Define the set $[n] := \{1, \dots, n\}$ for each $n\in \mathbb{P}$, and let $[0] := \emptyset$. Given ${\bf a} = (a_1, \dots, a_n)$, ${\bf b} = (b_1, \dots, b_n)$ in $\mathbb{Z}^n$, let $|{\bf a}| := a_1 + \dots + a_n$, and write ${\bf a} \leq {\bf b}$ if $a_i \leq b_i$ for all $i\in [n]$. For brevity, let ${\bf 1}_n$ denote the $n$-tuple $(1, \dots, 1)$.

Throughout this paper, $S:= \fieldk[x_1, \dots, x_n]$ is a standard $\mathbb{N}$-graded polynomial ring on $n$ variables over an infinite field $\fieldk$. All ideals, whether contained in $S$ or otherwise, are assumed to be homogeneous and proper, and for any minimal set of generators of a given ideal, we always assume the generators are homogeneous. Given $\mathcal{A} = \{p_1, \dots, p_r\}$ a collection of forms in $S$, write $\langle \mathcal{A} \rangle$ or $\langle p_1, \dots, p_r\rangle$ to mean the ideal of $S$ generated by $\mathcal{A}$. A {\it multicomplex} $M$ is a collection of monomials that is closed under divisibility (i.e. if $m\in M$ and $m^{\prime}$ divides $m$, then $m^{\prime} \in M$). For $\boldsymbol{\alpha} = (\alpha_1, \dots, \alpha_n) \in \mathbb{N}^n$, write ${\bf x}^{\boldsymbol{\alpha}}$ to mean the monomial ${\bf x}^{\boldsymbol{\alpha}} := x_1^{\alpha_1}\cdots x_n^{\alpha_n}$.

For any ring $R$ and any ideal $I\subset R$, let $\depth(I)$ denote the depth of $I$. We say $I$ is {\it perfect} if $\depth(I)$ equals the projective dimension of $R/I$. A {\it complete intersection ideal} (CI ideal) of $R$ is an ideal generated by an $R$-regular sequence of forms, and an {\it almost complete intersection ideal} of $R$ is a perfect ideal $I$ that is minimally generated by $\depth(I)+1$ elements. If $J = \langle f_1, \dots, f_r\rangle$ is a CI ideal such that $\deg(f_1) \leq \dots \leq \deg(f_r)$, then ${\bf e} := (\deg(f_1), \dots, \deg(f_r))\in \mathbb{P}^r$ is called the {\it type} of $J$ (as a CI ideal).

Given a finitely generated graded $S$-module $R = \bigoplus_{d\in \mathbb{N}} R_d$, let $H(R,-)$ denote its {\it Hilbert function}, i.e. $H(R,d) := \dim_{\fieldk} R_d$ for each $d\in \mathbb{N}$, and let $\soc(R)$ denote the socle of $R$. We say $R$ is {\it Gorenstein} if it is a Cohen-Macaulay module of Krull dimension $r$, and $\dim_{\fieldk} \soc(R/\langle f_1, \dots, f_r\rangle R) = 1$ for some $R$-regular sequence $f_1, \dots, f_r$. Note that complete intersection ideals are Gorenstein.

\defn Let $R$ be a Cohen-Macaulay ring, and let $I, I \subset R$ be ideals of height $r$. If there exists a CI ideal $J\subset R$ of height $r$ satisfying $J\subseteq I \cap I^{\prime}$, $I = J:I^{\prime}$ and $I^{\prime} = J:I$, then we say $I$ and $I^{\prime}$ are {\it (algebraically) directly linked} (by $J$), and we write $I \overset{J}{\sim} I^{\prime}$, or simply $I \sim I^{\prime}$ (if the ideal $J$ is not important). This binary relation $\sim$ is called {\it direct linkage}, and the ideal $J$ is called a {\it link}.

Direct linkage is symmetric, but not necessarily reflexive or transitive. Taking its transitive closure, we get an equivalence relation called {\it liaison} (or {\it linkage}). Equivalently, we have the following definition:

\defn Let $R$ be a Cohen-Macaulay ring, and let $I, I^{\prime}\subset R$ be ideals of height $r$. If there exists a finite sequence of ideals $I_0, I_1, \dots, I_s$ ($s\in \mathbb{P}$) of height $r$, such that $I_0 = I$, $I_s = I^{\prime}$, and $I_0 \sim I_1 \sim \dots \sim I_s$, then we call ``$I_0 \sim I_1 \sim \dots \sim I_s$'' a {\it sequence of links} from $I$ to $I^{\prime}$, we say $s$ is the {\it length} of this sequence, and we say $I$ and $I^{\prime}$ are {\it (algebraically) linked}. This binary relation (that $I$ and $I^{\prime}$ are linked) is an equivalence relation called {\it liaison} (or {\it linkage}), and each equivalence class is called a {\it liaison class}. An ideal that is in the liaison class\footnote{Any two CI ideals of the same height are linked, so for a fixed height, the liaison class of a CI ideal is unique; see, e.g. \cite{Schwartau1982:thesis}, for a proof.} of a CI ideal is called {\it licci}.

In our definition of liaison, we require that the links are CI ideals. This can be generalized by allowing links to be Gorenstein ideals, which yields the notion of {\it Gorenstein liaison}. The prefixes `CI-' and `G-' (for complete intersection and Gorenstein respectively) are usually attached to distinguish between the two definitions, e.g. CI-liaison, G-linked, etc.. Currently, Gorenstein liaison theory is an area of active research \cite{Gorla2008:GeneralizedGaetaThm,Hartshorne2007:GeneralizedDivisorsAndBiliaison,HartshorneSabadiniSchlesinger2008,KleppeMiglioreMiro-RoigNagelPeterson2001,MiglioreNagel2002:MonomialIdealsGLiaisonClassOfCI,MiglioreNagel2010:MinimalLinksResultGaeta}, and there are many open problems on whether results in CI-liaison theory can be extended analogously in G-liaison theory; see \cite{KleppeMiglioreMiro-RoigNagelPeterson2001}. In this paper, we only use CI-liaison theory and hence do not use the `CI-' prefix. For a good introduction to liaison theory, see \cite{book:MiglioreIntroLiaisonTheory}.

\remark The notion of liaison is (more commonly) defined for projective schemes as follows: Let $V_1, V_2$ be equidimensional subschemes of $\mathbb{P}^n_{\fieldk}$ of codimension $r$, and let $X$ be a complete intersection scheme of codimension $r$ containing both $V_1$ and $V_2$. If the defining ideals $I_{V_1}, I_{V_2}, I_X$ (of $V_1, V_2, X$ respectively) satisfy $I_X \subseteq I_{V_1} \cap I_{V_2}$, $I_X: I_{V_1} = I_{V_2}$ and $I_X: I_{V_2} = I_{V_1}$, then we say $V_1$ and $V_2$ are {\it (algebraically) directly linked}, and we write $V_1 \overset{X}{\sim} V_2$, or simply $V_1 \sim V_2$. The equivalence relation generated by (the transitive closure of) $\sim$ is called {\it liaison}, and we can define {\it link}, {\it liaison class}, etc. analogously.

\defn Let $R$ be a Cohen-Macaulay ring, and let $I, I^{\prime}, I^{\prime\prime} \subset R$ be ideals of height $r$. If $I$ minimally contains an $(a_1, \dots, a_r)$-regular sequence of forms, and $I$ and $I^{\prime}$ are directly linked by a CI ideal $J$ of type $(a_1, \dots, a_r)$, then we say $J$ is a {\it minimal} link. If $I$ and $I^{\prime\prime}$ are linked, and there exists a sequence of links from $I$ to $I^{\prime\prime}$ such that each link is minimal, then we say $I$ is {\it minimally linked to} $I^{\prime\prime}$. An ideal that is minimally linked to a CI ideal is called {\it minimally licci}.

Given linked ideals $I$ and $I^{\prime}$ of $S$, there are many possible sequences of links from $I$ to $I^{\prime}$ of varying lengths, and much work has been done on understanding licci ideals and their corresponding sequences of links. Gaeta~\cite{Gaeta1948} showed that every Cohen-Macaulay ideal of height two is minimally licci, while Peskine and Szpiro~\cite{PeskineSzpiro:Liaison} proved that an ideal of height two is licci if and only if it is Cohen-Macaulay. However, these results do not extend to ideals of height $\geq 3$. Not every Cohen-Macaulay ideal of height three is licci~\cite{HunekeUlrich1987:StructureOfLinkage}, and for every $r\geq 3$, there are licci ideals of height $r$ that are not minimally licci~\cite{HunekeMiglioreNagelUlrich2007}. Nevertheless, Watanabe~\cite{Watanabe1973:NoteOnGorenSteinRings} showed that every Gorenstein ideal of height three is licci, while Migliore and Nagel~\cite{MiglioreNagel2010:MinimalLinksResultGaeta} recently proved that every Gorenstein ideal of height three is minimally licci.

\section{Sequentially Bounded Licci Ideals}\label{sec:SeqBoundedLicci}
In this section, we weaken the notion of `minimally licci ideals' and define a subclass of licci ideals that we call `sequentially bounded'. In particular, minimally licci ideals are sequentially bounded licci ideals. As the main result of this section, we prove that the EGH conjecture holds for every sequentially bounded licci ideal that `admits a minimal first link', which we now define precisely:

\defn Let $I$ be a licci ideal of height $r$. Suppose there exists a sequence of links $I_0 \overset{J_1}{\sim} I_1 \overset{J_2}{\sim} \dots \overset{J_s}{\sim} I_s$ from $I_0 = I$ to a CI ideal $I_s$, such that $J_1, \dots, J_s$ (as CI ideals) have types ${\bf a}^{(1)}, \dots, {\bf a}^{(s)} \in \mathbb{P}^r$ respectively and satisfy ${\bf a}^{(1)} \geq \dots \geq {\bf a}^{(s)}$. Then $I$ is called a {\it sequentially bounded} licci ideal. Furthermore, if $J_1$ is a minimal link, then we say $I$ is a sequentially bounded licci ideal that {\it admits a minimal first link}.

\thm\label{thm:SeqBoundedLicci} Let $2\leq e_1 \leq \dots \leq e_n$ be integers. If $I\subsetneq S$ is a sequentially bounded licci ideal that admits a minimal first link and minimally contains an $(e_1, \dots, e_n)$-regular sequence of forms, then there exists a monomial ideal $J\subsetneq S$ containing $x_1^{e_1}, \dots, x_n^{e_n}$ such that $I$ and $J$ have the same Hilbert function.

Before we prove Theorem \ref{thm:SeqBoundedLicci}, we need the following useful theorem on Hilbert functions under liaison.
\thm[{\cite{DavisGeramitaOrecchia1985}}]\label{thm:HilbertFunctionsUnderLiaison}
Let $I,J$ be (homogeneous) ideals of $S$ such that $J\subseteq I$. If $S/J$ is an Artinian Gorenstein ring, and $s := \max\{t\in \mathbb{N}: H(S/J,t) \neq 0\}$, then $H(S/I,t) = H(S/J,t) - H(S/(J:I), s-t)$ for all $0\leq t\leq s$.

\begin{proof}[Proof of Theorem \ref{thm:SeqBoundedLicci}]
Let $I_0 \overset{J_1}{\sim} I_1 \overset{J_2}{\sim} \dots \overset{J_s}{\sim} I_s$ be a sequence of links from $I_0 = I$ to a CI ideal $I_s$, where $J_1$ is a minimal link. Let ${\bf a}^{(1)}, \dots, {\bf a}^{(s)}, {\bf a}^{(s+1)} \in \mathbb{P}^n$ denote the types of $J_1, \dots, J_s, I_s$ respectively (as CI ideals), and assume ${\bf a}^{(1)} \geq \dots \geq {\bf a}^{(s)}$. Note also that ${\bf a}^{(s)} \geq {\bf a}^{(s+1)}$, since $J_s$ and $I_s$ are CI ideals satisfying $J_s \subseteq I_s$. For each $i\in [s+1]$, write ${\bf a}^{(i)} = (a_{i,1}, \dots, a_{i,n})$, and let $\Gamma_i$ be the collection of all monomials in $S$ that divide ${\bf x}^{({\bf a}^{(i)} - {\bf 1}_n)}$. Also, define the collections of monomials $\widetilde{\Gamma}_{s+1}, \widetilde{\Gamma}_s, \dots, \widetilde{\Gamma}_1$ recursively as follows: Let $\widetilde{\Gamma}_{s+1} = \Gamma_{s+1}$, and for each $i\in [s]$, define
\begin{equation*}
\widetilde{\Gamma}_{s+1-i} := \Big\{q\in \Gamma_{s+1-i}: \frac{{\bf x}^{({\bf a}^{(s+1-i)} - {\bf 1}_n)}}{q} \not\in \widetilde{\Gamma}_{s+2-i}\Big\}.
\end{equation*}

\claimm $\widetilde{\Gamma}_i$ is a multicomplex, and $H(S/I_{i-1}, t) = \big|\big\{q\in \widetilde{\Gamma}_i: \deg(q) = t\big\}\big|$ for all $i\in [s+1]$, $t\in \mathbb{N}$.

\begin{proof}[Proof of Claim]
\renewcommand{\qedsymbol}{$\blacksquare$}
We shall prove both assertions of the claim simultaneously by induction on $s+2-i$ (for $i\in [s+1]$). The base case is trivial: $\widetilde{\Gamma}_{s+1}$ is clearly a multicomplex, and
\begin{equation*}
H(S/I_s,t) = H(S/\langle x_1^{a_{n+1,1}}, \dots, x_n^{a_{n+1,n}}\rangle,t) = \big|\big\{q\in \widetilde{\Gamma}_{s+1}: \deg(q) = t\big\}\big|
\end{equation*}
for all $t\in \mathbb{N}$. In particular, $\Gamma_{s+1} = \widetilde{\Gamma}_{s+1}$ forms a $\fieldk$-basis for $S/\langle x_1^{a_{n+1,1}}, \dots, x_n^{a_{n+1,n}}\rangle$.

For the induction step, let $m_{s+1-i} := {\bf x}^{({\bf a}^{(s+1-i)} - {\bf 1}_n)}$. Since ${\bf a}^{(1)} \geq \dots \geq {\bf a}^{(s+1)}$ implies $\widetilde{\Gamma}_{i+1} \subseteq \Gamma_{i+1} \subseteq \Gamma_i$ for all $i\in [s]$, it then follows from Theorem \ref{thm:HilbertFunctionsUnderLiaison} that $H(S/I_{s-i},t)$ equals
\begin{align*}
&\ H(S/J_{s+1-i},t) - H(S/I_{s+1-i}, |{\bf a}^{(s+1-i)}| -n-t)\\
=&\ \big|\big\{q\in \Gamma_{s+1-i}: \deg(q) = t\big\}\big| - \big|\big\{q^{\prime} \in \widetilde{\Gamma}_{s+2-i}: \deg(q^{\prime}) = |{\bf a}^{(s+1-i)}| -n-t\big\}\big|\\
=&\ \big|\big\{q\in \Gamma_{s+1-i}: \deg(q) = t\big\}\big| - \Big|\Big\{\frac{m_{s+1-i}}{q^{\prime}} \in \Gamma_{s+1-i}: q^{\prime} \in \widetilde{\Gamma}_{s+2-i}, \deg\Big(\frac{m_{s+1-i}}{q^{\prime}}\Big) = t\Big\}\Big|\\
=&\ \big|\big\{q\in \Gamma_{s+1-i}: \deg(q) = t\big\}\big| - \Big|\Big\{q\in \Gamma_{s+1-i}: \frac{m_{s+1-i}}{q} \in \widetilde{\Gamma}_{s+2-i}, \deg(q) = t\Big\}\Big|\\
=&\ \big|\big\{q\in \widetilde{\Gamma}_{s+1-i}: \deg(q) = t\big\}\big|
\end{align*}
for every $t\in \mathbb{N}$.

Next, we prove that $\widetilde{\Gamma}_{s+1-i}$ is a multicomplex. Choose an arbitrary $q\in \Gamma_{s+1-i}$ such that $q\not\in \widetilde{\Gamma}_{s+1-i}$, and suppose $qx_j \in \Gamma_{s+1-i}$ for some $j\in [n]$. Clearly $\Gamma_{s+1-i}$ is a multicomplex containing $\widetilde{\Gamma}_{s+1-i}$, hence to prove that $\widetilde{\Gamma}_{s+1-i}$ is a multicomplex, it suffices to show that $qx_j \not\in \widetilde{\Gamma}_{s+1-i}$. Now, $qx_j \in \Gamma_{s+1-i}$ implies $x_j^{(a_{(s+1-i),j}-1)}$ does not divide $q$, which means $x_j$ divides $\frac{m_{s+1-i}}{q}$, so since $\widetilde{\Gamma}_{s+2-i}$ is a multicomplex by induction hypothesis, we thus get $\frac{m_{s+1-i}}{qx_j} \in \widetilde{\Gamma}_{s+2-i}$ which yields $qx_j \not\in \widetilde{\Gamma}_{s+1-i}$.
\end{proof}

With the above claim, we now complete the proof of Theorem \ref{thm:SeqBoundedLicci}. Let $J\subset S$ be the ideal spanned by monomials in $S$ that are not contained in $\widetilde{\Gamma}_1$. Using the claim, $\widetilde{\Gamma}_1$ is a multicomplex, hence $\widetilde{\Gamma}_1$ forms a $\fieldk$-basis for $S/J$, and we get
\begin{equation*}
H(S/I,t) = H(S/I_0,t) = \big|\big\{q\in \widetilde{\Gamma}_1: \deg(q) = t\big\}\big| = H(S/J,t).
\end{equation*}
Finally, since $J_1$ is a minimal link, we have $(e_1, \dots, e_n) = {\bf a}^{(1)}$, so $\widetilde{\Gamma}_1 \subseteq \Gamma_1$ implies $J$ contains $x_1^{e_1}, \dots, x_n^{e_n}$.
\end{proof}

\cor\label{cor:minimallyLicci=>EGH} Let $2\leq e_1 \leq \dots \leq e_n$ be integers. If $I\subset S$ is a minimally licci ideal that minimally contains an $(e_1, \dots, e_n)$-regular sequence of forms, then there exists a monomial ideal in $S$ containing $x_1^{e_1}, \dots, x_n^{e_n}$ with the same Hilbert function as $I$.

\begin{proof}
Given any sequence $I_0 \overset{J_1}{\sim} I_1 \overset{J_2}{\sim} I_2$ of minimal links, where $J_1$ and $J_2$ have types ${\bf a}$ and ${\bf b}$ respectively, the definition of liaison yields $I_0 \overset{J_1}{\sim} I_1 \overset{J_1}{\sim} I_0$, hence the minimality of $J_2$ as a link implies ${\bf b} \leq {\bf a}$. Consequently, a minimally licci ideal is a sequentially bounded licci ideal that admits a minimal first link, and the assertion follows from Theorem \ref{thm:SeqBoundedLicci}.
\end{proof}

Theorem \ref{mainThm-Gorenstein} is an immediate consequence of Corollary \ref{cor:minimallyLicci=>EGH}, since every Gorenstein ideal of height three is miminally licci~\cite{MiglioreNagel2010:MinimalLinksResultGaeta}. In fact, Migliore and Nagel~\cite{MiglioreNagel2010:MinimalLinksResultGaeta} proved a stronger result: If $I$ is a Gorenstein ideal that is not a CI ideal, then linking $I$ minimally twice gives a Gorenstein ideal with two fewer generators than $I$. Their proof uses Buchsbaum-Eisenbud's structure theorem~\cite{BuchsbaumEisenbud1977}, which says that every Gorenstein ideal of height three is generated by the submaximal Pfaffians of an alternating matrix. Note that Watanabe~\cite{Watanabe1973:NoteOnGorenSteinRings} previously showed the minimal number of generators of every Gorenstein ideal of height three is odd.

\remark[The Two Variables Case] Every ideal $I\subset S$ containing a maximal $S$-regular sequence is Artinian and hence Cohen-Macaulay. Since Gaeta's theorem~\cite{Gaeta1948} says every Cohen-Macaulay ideal of height two is minimally licci, Corollary \ref{cor:minimallyLicci=>EGH} thus yields a different proof (cf. \cite{Richert2004:StudyLPPConjecture}, \cite[Remark 14]{CavigliaMaclagan2008:EGHConj}) that the EGH conjecture holds for $n=2$.

\section{Variations of the EGH Conjecture}
There are two common variations of the EGH conjecture.
\conjecture[{$\text{EGH}_{{\bf e},n}$}]\label{conj:non-minimal-degrees} Let $2\leq e_1 \leq \dots \leq e_n$ be integers, and let $f_1, \dots, f_n$ be an $S$-regular sequence of forms of degrees $e_1, \dots, e_n$ respectively. If $I\subsetneq S$ is a homogeneous ideal containing $f_1, \dots, f_n$, then there exists a homogeneous ideal $J\subsetneq S$ containing $x_1^{e_1}, \dots, x_n^{e_n}$, such that $I$ and $J$ have the same Hilbert function.

\conjecture[{$\text{EGH}_{n,{\bf e},r}$}]\label{conj:non-maximal-regularSeq} Let $r \in [n]$, ${\bf e} = (e_1, \dots, e_r) \in \mathbb{P}^r$ satisfy $2\leq e_1 \leq \dots \leq e_r$, and let $f_1, \dots, f_r$ be an $S$-regular sequence of forms of degrees $e_1, \dots, e_r$ respectively. If $I\subsetneq S$ is a homogeneous ideal containing $f_1, \dots, f_r$, then there exists a homogeneous ideal $J\subsetneq S$ containing $x_1^{e_1}, \dots, x_r^{e_r}$, such that $I$ and $J$ have the same Hilbert function.

Conjecture \ref{conj:non-minimal-degrees} is clearly equivalent to the EGH conjecture, while Conjecture \ref{conj:non-maximal-regularSeq} allows for non-maximal $S$-regular sequences, with the case $\text{EGH}_{n,{\bf e},n}$ being identical to $\text{EGH}_{{\bf e},n}$. Remarkably, Caviglia and Maclagan~\cite{CavigliaMaclagan2008:EGHConj} showed that Conjecture \ref{conj:non-maximal-regularSeq} is equivalent to the EGH conjecture. In particular, they showed that if $r\in [n]$, ${\bf e} = (e_1, \dots, e_r) \in \mathbb{P}^r$ satisfies $2\leq e_1 \leq \dots \leq e_r$, and $\text{EGH}_{{\bf e}^{\prime},n}$ holds for all ${\bf e}^{\prime} = (e_1^{\prime}, \dots, e_n^{\prime}) \in \mathbb{P}^n$ such that $2\leq e_1^{\prime} \leq \dots \leq e_n^{\prime}$ and $e_i^{\prime} = e_i$ for each $i\in [r]$, then $\text{EGH}_{n,{\bf e},r}$ holds; and conversely, if $\text{EGH}_{{\bf e},r}$ holds for some $r\in \mathbb{P}$ and some ${\bf e} = (e_1, \dots, e_r) \in \mathbb{P}^r$ satisfying $2\leq e_1 \leq \dots \leq e_r$, then $\text{EGH}_{n^{\prime},{\bf e},r}$ holds for all integers $n^{\prime} \geq r$.

By a modification of their proof, we show that $\text{EGH}_{n,{\bf e},r}$ holds for every sequentially bounded licci ideal that admits a minimal first link whenever $n \geq r$.

\lemma\label{lemma:non-maximalCase} Let $I_1, I_2 \subsetneq S$ be ideals of height $r<n$, and let $f_1, \dots, f_r,g$ be an $S$-regular sequence of forms. Define $J := \langle f_1, \dots, f_r\rangle$ and $J^{\prime} := \langle f_1, \dots, f_r, g\rangle$. If $I_1\overset{J}{\sim} I_2$, then
\begin{equation*}
\big((I_1: g^j) + \langle g \rangle\big) \overset{J^{\prime}}{\sim} \big((I_2: g^j) + \langle g \rangle\big)
\end{equation*}
for every $j\in \mathbb{N}$.

\begin{proof}
For convenience, write $I_i^{\prime} := \big((I_i: g^j) + \langle g \rangle\big)$ for each $i\in \{1,2\}$. The case $j=0$ is trivial, so assume $j\geq 1$. Observe that the identity $I_2 = J:I_1$ yields
\begin{align*}
I_2^{\prime} &= \frac{1}{g^j}\Big((J:I_1) \cap \langle g^j\rangle\Big) + \langle g\rangle = \Big\{\frac{s^{\prime}}{g^j} \in S: s^{\prime} \in \langle g^j\rangle, s^{\prime}I_1 \subseteq J\Big\} + \langle g \rangle\\
&= \{s\in S: sg^jI_1\subseteq J\} + \langle g\rangle,
\end{align*}
while the identity $I_1 = J:(J:I_1)$ yields
\begin{align*}
J^{\prime}: I_1^{\prime} &= \Big\{s\in S: s\cdot \frac{1}{g^j}(I_1\cap \langle g^j\rangle) \subseteq J\Big\} + \langle g\rangle\\
&= \Big\{s\in S: s\cdot \frac{1}{g^j}((J:(J:I_1))\cap \langle g^j\rangle) \subseteq J\Big\} + \langle g\rangle\\
&= \big\{s\in S: ss^{\prime} \in J\text{ for every }s^{\prime}\in S^{\prime}\big\} + \langle g\rangle,
\end{align*}
where $S^{\prime} := \{s^{\prime} \in S: s^{\prime}tg^j \in J\text{ for all }t\in S \text{ such that }tI_1\subseteq J\}$.

A routine check gives $J^{\prime}: I_1^{\prime} \subseteq I_2^{\prime}$. To show the reverse inclusion $J^{\prime}: I_1^{\prime} \supseteq I_2^{\prime}$, note that $sp\in J$ for every non-zero $s\in S$ and non-zero $p\in I_1$ satisfying $sg^jp\in J$, since otherwise $g^j$ would be a zero-divisor of $S/J$, which contradicts the assumption that $f_1, \dots, f_r,g$ is an $S$-regular sequence. By the same argument, every $s^{\prime} \in S^{\prime}$ satisfies $s^{\prime}t$ for all $t\in S$ such that $tI_1\subseteq J$. Consequently, each $s \in I_2^{\prime}$ satisfies $sI_1\subseteq J$ and hence satisfies $ss^{\prime} \in J$ for all $s^{\prime} \in S^{\prime}$, which gives the reverse inclusion, so we conclude $J^{\prime}: I_1^{\prime} = I_2^{\prime}$. A symmetric argument yields $J^{\prime}: I_2^{\prime} = I_1^{\prime}$.

Finally, since $f_1, \dots, f_r, g$ is an $S$-regular sequence, we get $J\cap \langle g^j\rangle = \langle g^jf_1, \dots, g^jf_r\rangle$. It then follows from $J\subseteq I_1 \cap I_2$ that $J^{\prime} = \frac{1}{g^j}(J \cap \langle g^j\rangle) + \langle g\rangle \subseteq I_1^{\prime} \cap I_2^{\prime}$, therefore $I_1^{\prime} \overset{J^{\prime}}{\sim} I_2^{\prime}$.
\end{proof}

\prop Let $r\in [n]$, let $2\leq e_1 \leq \dots \leq e_r$ be integers, and let $I\subsetneq S$ be an ideal (of height $r$) that minimally contains an $(e_1, \dots, e_r)$-regular sequence of forms. If $I\subsetneq S$ is a sequentially bounded licci ideal that admits a minimal first link (for example, $I$ could be a minimally licci ideal), then there exists a monomial ideal $J\subsetneq S$ containing $x_1^{e_1}, \dots, x_r^{e_r}$ such that $I$ and $J$ have the same Hilbert function.

\begin{proof}
We follow Caviglia-Maclagan's proof and similarly prove our proposition by induction on $n-r\in \mathbb{N}$, where the base case $n-r=0$ is equivalent to Theorem \ref{thm:SeqBoundedLicci}. Assume $r <n$, fix the integers $2\leq e_1 \leq \dots \leq e_r$, and let $I = I_0 \overset{J_1}{\sim} I_1 \overset{J_2}{\sim} \dots \overset{J_s}{\sim} I_s$ be a sequence of links from $I_0$ to a CI ideal $I_s$, such that $J_1, \dots, J_s, I_s$ (as CI ideals) have types ${\bf a}^{(1)}, \dots, {\bf a}^{(s)}, {\bf a}^{(s+1)} \in \mathbb{P}^r$ respectively and satisfy ${\bf a}^{(1)} \geq \dots \geq {\bf a}^{(s+1)}$. In particular, $J_s$ and $I_s$ are CI ideals satisfying $J_s \subseteq I_s$, so the last inequality ${\bf a}^{(s)} \geq {\bf a}^{(s+1)}$ is guaranteed. Also, assume $J_1$ is a minimal link and ${\bf a}^{(1)} = (e_1, \dots, e_r)$. Since $\fieldk$ is infinite and $r<n$, we can choose some linear form $g$ that is a non-zero-divisor on each of $S/J_1, \dots, S/J_s, S/I_s$. Consequently, for each $i\in [s]$, we can define the CI ideal $J_i^{\prime} := J_i + \langle g\rangle$ of type $(1,{\bf a}^{(i)}) \in \mathbb{P}^{r+1}$.

Next, let $N\in \mathbb{P}$ be sufficiently large so that $(I:g^{\infty}) = (I:g^N)$. For each $i\in \{0,1,\dots, s\}$ and $j\in \{0,1,\dots, N\}$, define the ideal $I_i^{(j)} := (I_i: g^j) + \langle g\rangle$. By Lemma \ref{lemma:non-maximalCase}, we get the sequence of links
\begin{equation}\label{eqn:seqLinks}
I_0^{(j)} \overset{J_1^{\prime}}{\sim} I_1^{(j)} \overset{J_2^{\prime}}{\sim} \dots \overset{J_s^{\prime}}{\sim} I_s^{(j)}
\end{equation}
for every $j\in \{0,1,\dots, N\}$, and we observe that each $I_s^{(j)}$ is a CI ideal.

For each $i\in [s]$, write $J_i = \langle f_1^{(i)}, \dots, f_r^{(i)} \rangle$ so that ${\bf a}^{(i)} = (\deg(f_1^{(i)}), \dots, \deg(f_r^{(i)}))$. Also, write $I_s = \langle f_1^{(s+1)}, \dots, f_r^{(s+1)}\rangle$ so that ${\bf a}^{(s+1)} = (\deg(f_1^{(s+1)}), \dots, \deg(f_r^{(s+1)}))$. The quotient ring $R := S/\langle g\rangle$ is isomorphic to a polynomial ring on $n-1$ variables, and by the natural quotient map $\pi: S\to R$, the $S$-regular sequence $f_1^{(i)}, \dots, f_r^{(i)}$ descends to an $R$-regular sequence $\pi(f_1^{(i)}), \dots, \pi(f_r^{(i)})$ for each $i\in [s+1]$. In particular, $J_i^{\prime\prime} := \pi(J_i^{\prime})$ is a CI ideal in $R$ for all $i\in [s]$. By applying $\pi$ to the ideals in \eqref{eqn:seqLinks}, we get the sequence of links $\pi(I_0^{(j)}) \overset{J_1^{\prime\prime}}{\sim} \pi(I_1^{(j)}) \overset{J_2^{\prime\prime}}{\sim} \dots \overset{J_s^{\prime\prime}}{\sim} \pi(I_s^{(j)})$ for every $j\in \{0,1,\dots, N\}$.

Now $J_1^{\prime\prime}, \dots, J_s^{\prime\prime}$ (as CI ideals) have types ${\bf a}^{(1)}, \dots, {\bf a}^{(s)}$ respectively, and for every $j\in \{0,1,\dots, N\}$, the ideal $\pi(I_s^{(j)})$ is a CI ideal of type ${\bf a}^{(s+1)}$, thus $\pi(I_0^{(j)})$ is a sequentially bounded licci ideal in $R$ that admits a minimal first link. Using $R \cong \fieldk[x_1,\dots, x_{n-1}]$, the induction hypothesis then says that for every $j\in \{0,1,\dots, N\}$, there exists a monomial ideal in $\fieldk[x_1,\dots, x_{n-1}]$ containing $x_1^{e_1}, \dots, x_r^{e_r}$ with the same Hilbert function as $\pi(I_0^{(j)})$. Let $M_j$ be the lex-plus-powers ideal in $\fieldk[x_1, \dots, x_{n-1}]$ containing $x_1^{e_1}, \dots, x_r^{e_r}$ with this same Hilbert function (which exists by Clements-Lindstr\"{o}m theorem~\cite{ClementsLindstrom1969}), and let $K_j \subseteq \fieldk[x_1, \dots, x_n]$ be the set of monomials $K_j := \{mx_n^j: m\in M_j\}$.

Finally, consider the ideal $K$ generated by the monomials in $\bigcup_{j=0}^N K_j$. As shown by Caviglia and Maclagan~\cite[Proposition 10]{CavigliaMaclagan2008:EGHConj}, the ideal $K$ contains $x_1^{e_1}, \dots, x_r^{e_r}$ and has the same Hilbert function as $I$, and their proof holds verbatim.
\end{proof}

\cor Let $n\geq 3$, and let $2\leq e_1 \leq e_2 \leq e_3$ be integers. If $I\subsetneq S$ is a homogeneous Gorenstein ideal (of height three) that minimally contains an $(e_1, e_2, e_3)$-regular sequence of forms, then there exists a monomial ideal $J \subsetneq S$ containing $x_1^{e_1}, x_2^{e_2}, x_3^{e_3}$, such that $I$ and $J$ have the same Hilbert function.

\section*{Acknowledgements}
The author thanks David Eisenbud, Irena Peeva, Juan Migliore, and Edward Swartz for helpful comments.

\bibliographystyle{plain}
\bibliography{References}

\begin{thebibliography}{10}

\bibitem{Abedelfatah2012:EGHConjPreprint}
Abed Abedelfatah.
\newblock On the {E}isenbud-{G}reen-{H}arris conjecture.
\newblock Preprint arXiv:1212.2653v1 [math.AC], December 2012.

\bibitem{BuchsbaumEisenbud1977}
David~A. Buchsbaum and David Eisenbud.
\newblock Algebra structures for finite free resolutions, and some structure
  theorems for ideals of codimension {$3$}.
\newblock {\em Amer. J. Math.}, 99(3):447--485, 1977.

\bibitem{CavigliaMaclagan2008:EGHConj}
Giulio Caviglia and Diane Maclagan.
\newblock Some cases of the {E}isenbud-{G}reen-{H}arris conjecture.
\newblock {\em Math. Res. Lett.}, 15(3):427--433, 2008.

\bibitem{Chen2011:thesis}
Ri-Xiang Chen.
\newblock {\em Hilbert functions and free resolutions}.
\newblock PhD thesis, Cornell University, 2011.

\bibitem{ClementsLindstrom1969}
G.~F. Clements and B.~Lindstr{\"o}m.
\newblock A generalization of a combinatorial theorem of {M}acaulay.
\newblock {\em J. Combinatorial Theory}, 7:230--238, 1969.

\bibitem{Cooper:EGHPreprintCase3}
Susan~M. Cooper.
\newblock The {E}isenbud-{G}reen-{H}arris conjecture for ideals of points.
\newblock Preprint.

\bibitem{DavisGeramitaOrecchia1985}
E.~D. Davis, A.~V. Geramita, and F.~Orecchia.
\newblock Gorenstein algebras and the {C}ayley-{B}acharach theorem.
\newblock {\em Proc. Amer. Math. Soc.}, 93(4):593--597, 1985.

\bibitem{EisenbudGreenHarris1993:EGHConjecture}
David Eisenbud, Mark Green, and Joe Harris.
\newblock Higher {C}astelnuovo theory.
\newblock {\em Ast\'erisque}, (218):187--202, 1993.
\newblock Journ{\'e}es de G{\'e}om{\'e}trie Alg{\'e}brique d'Orsay (Orsay,
  1992).

\bibitem{Francisco2004:AlmostCIandLPPConj}
Christopher~A. Francisco.
\newblock Almost complete intersections and the lex-plus-powers conjecture.
\newblock {\em J. Algebra}, 276(2):737--760, 2004.

\bibitem{FranciscoRichert2007:LPPIdeals}
Christopher~A. Francisco and Benjamin~P. Richert.
\newblock Lex-plus-powers ideals.
\newblock In {\em Syzygies and {H}ilbert functions}, volume 254 of {\em Lect.
  Notes Pure Appl. Math.}, pages 113--144. Chapman \& Hall/CRC, Boca Raton, FL,
  2007.

\bibitem{Gaeta1948}
Federico Gaeta.
\newblock Sulle curve sghembe algebriche di residuale finito.
\newblock {\em Ann. Mat. Pura Appl. (4)}, 27:177--241, 1948.

\bibitem{Gasharov1999:HilbertFunctionsHomogGenericForms}
Vesselin Gasharov.
\newblock Hilbert functions and homogeneous generic forms. {II}.
\newblock {\em Compositio Math.}, 116(2):167--172, 1999.

\bibitem{Gorla2008:GeneralizedGaetaThm}
Elisa Gorla.
\newblock A generalized {G}aeta's theorem.
\newblock {\em Compos. Math.}, 144(3):689--704, 2008.

\bibitem{Hartshorne2007:GeneralizedDivisorsAndBiliaison}
Robin Hartshorne.
\newblock Generalized divisors and biliaison.
\newblock {\em Illinois J. Math.}, 51(1):83--98 (electronic), 2007.

\bibitem{HartshorneSabadiniSchlesinger2008}
Robin Hartshorne, Irene Sabadini, and Enrico Schlesinger.
\newblock Codimension 3 arithmetically {G}orenstein subschemes of projective
  {$N$}-space.
\newblock {\em Ann. Inst. Fourier (Grenoble)}, 58(6):2037--2073, 2008.

\bibitem{HerzogPopescu1998:HilbertFunctionsGenericForms}
J{\"u}rgen Herzog and Dorin Popescu.
\newblock Hilbert functions and generic forms.
\newblock {\em Compositio Math.}, 113(1):1--22, 1998.

\bibitem{HunekeMiglioreNagelUlrich2007}
Craig Huneke, Juan Migliore, Uwe Nagel, and Bernd Ulrich.
\newblock Minimal homogeneous liaison and licci ideals.
\newblock In {\em Algebra, geometry and their interactions}, volume 448 of {\em
  Contemp. Math.}, pages 129--139. Amer. Math. Soc., Providence, RI, 2007.

\bibitem{HunekeUlrich1987:StructureOfLinkage}
Craig Huneke and Bernd Ulrich.
\newblock The structure of linkage.
\newblock {\em Ann. of Math. (2)}, 126(2):277--334, 1987.

\bibitem{KleppeMiglioreMiro-RoigNagelPeterson2001}
Jan~O. Kleppe, Juan~C. Migliore, Rosa Mir{\'o}-Roig, Uwe Nagel, and Chris
  Peterson.
\newblock Gorenstein liaison, complete intersection liaison invariants and
  unobstructedness.
\newblock {\em Mem. Amer. Math. Soc.}, 154(732):viii+116, 2001.

\bibitem{Macaulay1927}
F.~S. Macaulay.
\newblock Some properties of enumeration in the theory of modular systems.
\newblock {\em Proc. London Math. Soc.}, 26(1):531--555, 1927.

\bibitem{MiglioreNagel2002:MonomialIdealsGLiaisonClassOfCI}
J.~Migliore and U.~Nagel.
\newblock Monomial ideals and the {G}orenstein liaison class of a complete
  intersection.
\newblock {\em Compositio Math.}, 133(1):25--36, 2002.

\bibitem{MiglioreNagel2010:MinimalLinksResultGaeta}
Juan Migliore and Uwe Nagel.
\newblock Minimal links and a result of {G}aeta.
\newblock In {\em Liaison, {S}chottky problem and invariant theory}, volume 280
  of {\em Progr. Math.}, pages 103--132. Birkh\"auser Verlag, Basel, 2010.

\bibitem{book:MiglioreIntroLiaisonTheory}
Juan~C. Migliore.
\newblock {\em Introduction to liaison theory and deficiency modules}, volume
  165 of {\em Progress in Mathematics}.
\newblock Birkh\"auser Boston Inc., Boston, MA, 1998.

\bibitem{PeskineSzpiro:Liaison}
C.~Peskine and L.~Szpiro.
\newblock Liaison des vari\'et\'es alg\'ebriques. {I}.
\newblock {\em Invent. Math.}, 26:271--302, 1974.

\bibitem{Richert2004:StudyLPPConjecture}
Benjamin~P. Richert.
\newblock A study of the lex plus powers conjecture.
\newblock {\em J. Pure Appl. Algebra}, 186(2):169--183, 2004.

\bibitem{Schwartau1982:thesis}
Philip~William Schwartau.
\newblock {\em Liaison addition and monomial ideals}.
\newblock ProQuest LLC, Ann Arbor, MI, 1982.
\newblock Thesis (Ph.D.)--Brandeis University.

\bibitem{Watanabe1973:NoteOnGorenSteinRings}
Junzo Watanabe.
\newblock A note on {G}orenstein rings of embedding codimension three.
\newblock {\em Nagoya Math. J.}, 50:227--232, 1973.

\end{thebibliography}

\end{document}